\documentclass[12pt]{article}
\usepackage{amsmath}
\usepackage{amsthm}
\usepackage{graphicx}
\usepackage{subcaption}
\usepackage{lineno}
\usepackage{verbatim}
\usepackage{epsfig}
\usepackage{hyperref}
\usepackage{float}
\usepackage{color}
\usepackage{fullpage}
\usepackage{amsfonts}
\usepackage{amscd}
\usepackage{amssymb}
\usepackage{graphics}
\usepackage{amsmath, bm}
\usepackage{amsbsy}
\usepackage[english]{babel}
\usepackage{bbm}
\usepackage{yfonts}
\usepackage{mathrsfs}
\usepackage{color}
\usepackage{enumitem}
\usepackage{cleveref}
\usepackage{authblk}
\usepackage{natbib}
\DeclareFontFamily{OML}{rsfs}{\skewchar\font'177}
\DeclareFontShape{OML}{rsfs}{m}{n}{ <5> <6> rsfs5 <7> <8> <9>
rsfs7 <10> <10.95> <12> <14.4> <17.28> <20.74> <24.88> rsfs10 }{}
\DeclareMathAlphabet{\mathfs}{OML}{rsfs}{m}{n}
\newcommand{\BE}{{\mathbb{E}}}

\newcommand{\BP}{{\mathbb{P}}}

\newcommand{\BR}{{\mathbb{R}}}

\newcommand{\CN}{{\mathcal{N}}}

\newcommand{\bae}{\begin{equation}\begin{aligned}}
\newcommand{\eae}{\end{aligned}\end{equation}}

\newtheorem{theorem}{Theorem}[section]

\theoremstyle{definition}

\begin{document}
\numberwithin{equation}{section} \numberwithin{figure}{section}
\title{Reconciling common source, specific source, feature based and score based likelihood ratios}%; it is all a matter of information}
\author{Aafko Boonstra\footnote{Department of Mathematics, Vrije Universiteit Amsterdam, De Boelelaan 1111, 1081 HV, Amsterdam, The Netherlands, a.j.w.boonstra@vu.nl.}, \, Ronald Meester\footnote{Department of Mathematics, Vrije Universiteit Amsterdam, De Boelelaan 1111, 1081 HV, Amsterdam, The Netherlands, r.w.j.meester@vu.nl.} \, and Klaas Slooten\footnote{Netherlands Forensic Institute, Department of Human Biological traces, Laan van Ypenburg 6, 2497 GB, The Hague, The Netherlands; Department of Mathematics, Vrije Universiteit Amsterdam, De Boelelaan 1111, 1081 HV, Amsterdam, The Netherlands; k.slooten@nfi.nl.}}
\maketitle

\begin{abstract}
We show that the incorporation of any new piece of information allows for improved decision making in the sense that the expected costs of an optimal decision decrease (or, in boundary cases where no or not enough new information is incorporated, stays the same) whenever this is done by the appropriate update of the probabilities of the hypotheses. 
Versions of this result have been stated before. However, previous proofs rely on auxiliary constructions with proper scoring rules. We, instead, offer a direct and completely general proof by considering elementary properties of likelihood ratios only. We apply our results to make a contribution to the debates about the use of score based/feature based and common/specific source likelihood ratios. In the literature these are often presented as different ``LR-systems''. We argue that the difference between these is simply a matter which  information is processed. There is no therefore no such thing as different ``LR-systems'', there are only differences in the processed information. In particular, despite claims to the contrary, scores can very well be used in forensic practice and we illustrate this with an extensive example in DNA kinship context.
\end{abstract}

%We contribute to the debates about the use of score based/feature based and common/specific source likelihood ratios. In the literature these are often presented as different “LR-systems”, processing the same information for different pairs of hypotheses.
%Differences between these LRs have been used to argue against the application of all score-based and common source LRs.
%We argue for the opposite point of view: these LRs are for the same hypotheses, but processing different information. Thus, there is no such thing as different “LR-systems”. The differences in LR are due to differences in the available information, but there is only one statistical model. In particular, despite claims to the contrary, we see no reason to categorically reject score-based or common-source LRs. We show that in fact, their use is widespread and (we believe rightly so) undisputed in the DNA context and illustrate this with extensive examples.

\section{Introduction, context and background}
\label{intro}
In forensics as well as in many other frameworks, one often needs to decide on an action to take in the face of uncertainty. Typically, one has several possible ground truths, one of which is actually true, but it is not known with certainty which one. If we, somehow, would know which one, then we would take some decision or choose some consequential action. One may think of, for example, the decision to give a patient a certain treatment, to decide in favor or against paternity of an alleged father, to convict (or not) a suspect, etc. We suppose that the different ground truths that we, a priori, deem possible are formulated as hypotheses $H_1, \dots, H_n$. In the absence of further knowledge, we let $\BP(H_i)$ stand for the (subjective) probability, as assigned by the observer tasked with the decision making, that $H_i$ is true. We assume here that the $H_i$ are mutually exclusive and exhaustive so that precisely one of them is true. Subjective probabilities, when rationally expressed (cf. \cite{meesterslootenbook} for the meaning of this), follow the Kolmogorov axioms so we can use standard probability theory.

In order to facilitate the decision making, we assume some data $e$ become known. Ideally, the probabilities $\BP(H_i)$ will be updated to $\BP(H_i \mid e)$, by calculation of all likelihoods $\BP(e\mid H_i)$ and application of Bayes rule \cite{enfsi, aitkentaroni}. In forensics an often encountered question is whether two measurements, the origin of at least one of them being unknown, are actually measurements from one and the same source. For example, the question could be whether two fingerprints were left by the same unknown finger \cite{abraham}, whether a particle of glass comes from a reference window pane \cite{prguide}, whether a certain bullet was fired by a certain gun \cite{matzen, vanderplas}, whether a trace fingerprint was left by a certain finger of a person we have reference prints of \cite{alberink2014}, etc. Similar questions arise in many other forensic disciplines. 

We denote the two measurements by $(e_x, e_y)$. These are often obtained from `raw' data $(d_x, d_y)$, by application of some protocol that cleans, discretizes, projects, or otherwise simplifies the data. A very basic example would correspond to rounding or binning of a measured continuous variable. Another example would consist of the application of detection thresholds, stutter filter, artefact pruning etc for a DNA profile \cite{rakay, westen}. In other words, we assume here that $(d_x, d_y)$ are the raw data, and that $(e_x, e_y)$ are the data after standard data cleaning procedures. We want to assess the impact of the pair $(e_x, e_y)$ on our probabilistic assessment of the truth of the $H_i$.

Suppose we let $(S_x, S_y)$ be the (physical) sources of $e=(e_x, e_y)$, in whatever population we have in mind that these measurements are obtained from. A standard question is whether $H_1: S_x =S_y$ or $H_2: S_x \neq S_y$ is true. Summarizing into $I$ any prior information about $S_x, S_y$, if we have access to $\BP(e \mid H_1, I)$ and $\BP(e \mid H_2, I)$, it suffices to compute the likelihood ratio $LR_{H_1,H_2}(e)=\BP(e\mid H_1, I)/\BP(e \mid H_2, I)$ to obtain the posterior odds $\BP(H_1 \mid e, I)/\BP(H_2 \mid e, I)=LR_{H_1,H_2}(e)\BP(H_1 \mid I)/\BP(H_2 \mid I)$. In practice, this is not straightforward.

Somewhat as an aside, we mention that we do not reserve the term `likelihood ratio' for hypotheses where all statistical parameters have a fixed value, but also use it for the cases where they are specified via any probability distribution. In the general statistical literature, the latter is often called `Bayes factor' and the likelihood ratio is sometimes considered as the corresponding function of the parameters. In forensics, it is more common to use the term LR for $P(E\mid H_1)/P(E\mid H_2)$ irrespective of whether all parameters are fixed or not. To us, this is logical, as the value of evidence as obtained by an observer will depend on the observer's knowledge and understanding of the parameters, and conditional on that, not on the parameters themselves. The case where this knowledge leads to a point distribution is in that philosophy an unremarkable special case of uncertainty.  %Regardless of the extent of our knowledge, we must assign the likelihoods $P(E \mid H_i)$ in order to assign the evidential value. 

In general, when setting up a statistical model for the measurements we need to distinguish between the distribution of data that a known fixed source can yield, and the variation that exists between sources. To make this more precise, we consider measurements on a source $S$ to be, for the data generation process, possibly inherently random but probabilistically determined by some parameters $\theta_S$ pertaining to the source itself, and a measurement model predicting the observations. For example, in the DNA context a source could be an individual, and the parameters $\theta_S$ of that individual could be his or her DNA profile. The measurements of that profile would deterministically obtain that profile in case we model a reference sample, or have some random variation in the form of alleles randomly dropping out of or into the measured profile, if we model a trace profile left by that individual. For glass comparison purposes, a source would correspond to a glass object, and its parameters $\theta_S$ the true values for that glass object of what it is we measure (e.g., its refractive index, or elemental composition for the measured elements). Note that there is a distinction between having an identified source $S$ (i.e., the physical identification of $S$) and knowing its parameters $\theta_S$. For the probability distribution of measurements on $S$, given $\theta_S$ the identity of $S$ is irrelevant. 
For example, if we know someone's DNA profile we consider this a known source, regardless of whether or not we know the identity of that person. 

Typically, at least one of the hypotheses also involves unknown sources, so that we are forced to consider the distribution of the source parameters, usually (but not necessarily) by considering these distributions as representative for some population of sources. For example, in the DNA context we will consider the distribution of DNA profiles, e.g., by measuring allele frequencies and assuming Hardy-Weinberg equilibrium. This would give a description of (the DNA profiles of) unknown individuals unrelated to any known individuals.
In the glass context we would need to consider a population of sources (i.e., glass objects) that are candidate sources for the unknown sources. The evaluation of a likelihood $\BP(e \mid H)$ might require the integration over the unknown sources.

Recall that $(e_x, e_y)$ are the data that we have of sources $S_x, S_y$. It may be that $e_x=\theta_{S_x}$, i.e., that $S_x$ is a known source (we know the parameters of $S_x$). It may also be that we do not know $\theta_{S_x}$ but that $e_x$ represents some measurement that is informative for $\theta_{S_x}$. Regardless, we denote by $e_x$ the data pertaining to $S_x$. To simplify the evaluation of the data, the pair $(e_x, e_y)$ may be reduced further to, say, $g(e_x, e_y)$. This function may take the form $(h(e_x), h(e_y))$, e.g., for DNA traces we can let $h(e_x)$ correspond to omitting Y-chromosomal information and/or peak height information. In such a case, we might have also called $e_x$ the raw data, and $h(e_x)$ the actual data that we process. If the function $g$ computes a real number, it is customary to call this a score function. A score function is typically constructed so as to measure similarity between $e_x$ and $e_y$. For example, given two DNA profiles $(e_x, e_y)$ the score function could count the number of alleles the two profiles have in common, or even simply be $\delta_{e_x,e_y}$ measuring whether the profiles are identical or not.

In the forensic literature, different nomenclatures are used to distinguish between different likelihood ratio calculations. For example, a `common source' LR is one where the hypotheses state that $e_x$ and $e_y$ have the same unknown source versus two different unknown sources. A `specific source' LR is one where a source and its parameters are known, say $e_x=\theta_{S_x}$  \cite{ommen, saunders}. It is, of course, also conceivable that we have two known sources $S_1, S_2$, with the question being, e.g., whether $e_x, e_y$ are both from $S_1$ or both from $S_2$. In practice this is a less encountered situation. 

Another distinction is whether the LR is for the data $(e_x,e_y)$ or for $g(e_x,e_y)$. The first case, fully modeling $(e_x, e_y)$, is called a `feature based' approach. If, instead, we calculate a score $g(e_x, e_y)$ and then evaluate the evidential value of that score, this is called a `score based' approach. Clearly, a score usually carries less information. There are, of course, many intermediate situations possible where some $g(e_x, e_y)$ is evaluated as evidence, that carries less information than $(e_x, e_y)$ but that is not a real valued function, e.g., simply omitting some part of the data. Here, we will call these score based approaches as well.

Finally, in addition to the score, we may also consider only one of the pair $(e_x, e_y)$, e.g., $(e_x, g(e_x, e_y))$. If $e_x$ stands for the parameters of the known source $S_x$, and the score function is some similarity score, we will then evaluate how likely it is to find the observed similarity between $S_x$ and measurements $e_y.$ Here, we compare the hypothesis that $e_y$ is a measurement on $S_x$ to the hypothesis that $e_y$ is a measurement on some unknown source. If $S_x$ is a very typical source, similarity with $e_y$ will be less indicative for $e_y$ coming from $S_x$ compared to the case where $S_x$ has rare parameters. For both hypotheses, $S_x$ is the undisputed source of $e_x,$ so that the likelihood ratio based on $(e_x, g(e_x, e_y))$ can be obtained by computing $\BP(g(e_x, e_y) \mid H_i, e_x),$ i.e. by conditioning on $e_x$. These likelihood ratios are therefore said to be `anchored' on source $S_x$. 

All these different LRs correspond to different information positions that one might have for investigating whether $S_x=S_y$ or not.
In the literature, we see a tendency to view these different LRs as distinct in the sense that they all address their own set of hypotheses (\cite{vergeer} and references therein). This point of view implies that one needs to select the most relevant pair of hypotheses for one's data. Here, we argue instead that all these LRs address the same hypotheses, namely whether $S_x=S_y$ or not, but that the data that are available (or the data that are actually taken into account), differ. E.g., processing a score $g(e_x, e_y)$ will often lead to a different LR compared to $(e_x, e_y)$ and the latter approach certainly can be said to be a better informed probabilistic update on the hypotheses.  
That is not to say that the update with $g(e_x, e_y)$ is incorrect, only that it is sub-optimal from an information perspective having only incorporated part of the available information.

Keeping the goal in mind that one wants to take some action depending on the probabilistic assessment of the truth of the competing hypotheses, given any possible ground truth $H_i$, some actions are more desirable than others. If we assume that the desirability of actions given hypotheses can be assigned costs, by expressing costs $c_{i,j}$ as the cost of action $i$ when, in reality, hypothesis $j$ is true, then we can compare expected costs given different decision strategies.

In this article, we contribute to the discussion about the use of score based methods, and about the difference between the various scenarios in terms of common or specific source. We argue that one should approach these issues from an information-theoretic perspective.

We first prove a very general inequality (cf. Theorem \ref{bdimp}) which tells us that incorporating more evidence is on average better than not incorporating it. 
Previous proofs of various versions of this statement relied on auxiliary constructions with proper scoring rules \cite{degroot, dawid2007, brummer2006, raftery, brummer, ferrer, vergeer}. We, instead, offer a direct and completely general proof by considering elementary properties of likelihood ratios only. We believe that our approach will contribute to the understanding of Bayesian decision making for forensic and possibly other scientists. 

We apply this to score/feature based LRs and to the common/specific source discussion. Although score based methods are widely used in forensic science \cite{leegwater2017, hepler, annabel, gonzalez}, Neumann and Ausdemore \cite{neumann} argue that score based methods should not be used. They warn against the use of score based methods by constructing scenarios in which the score gives misleading information. Our arguments show that they are mistaken at this point.

In the common/specific source debate, we apply our results to show that there is no reason to see these as different ``LR-systems'': their difference is only a matter of which information is processed, whereas the underlying statistical models used are exactly the same. 

Summarizing, we offer a unified framework within which we can interpret and understand all issues arising in the score/feature and common/specific source debate.

\section{Bayes Decisions improve with more information}
\label{BD}
We start with the rather general situation in which we must decide on some action to take, such that our choice of action would depend on which one of a set of hypotheses is true, but we are not certain about which one is actually true. That is, we have some probability distribution  on a mutually exclusive and exhaustive set of hypotheses $H_1, \ldots, H_n$, and we must choose some action $A \in \{A_1, \dots, A_m\}$. Given the truth of a $H_i$, we consider some of the actions more desirable than others, which we express in terms of a \textit{cost function}; $c_{ij} \in \BR$ representing the cost of choosing action $A_i$ if in fact $H_j$ is true.  Cost functions are a widely used tool, including the field of forensic science \cite{vanlierop2024}. A special case arises when  $m=n$ while $A_i$ is identified with the decision that $H_i$ is true (a so called `hard decision'), but we do not require this, nor do we suppose anything about the $c_{ij}$.

A rational approach is to choose the action that we expect to have lowest cost in light of our probability distribution $\pi=(\pi_1,\dots,\pi_n) :=(\BP(H_1),\dots, \BP(H_n))$. Given this distribution, the expected cost of action $A_i$ is 
\begin{equation}
\label{expc}
c_\pi(A_i)=\sum_{j=1}^n c_{ij}\pi_j.
\end{equation}

Any action with lowest cost will be referred to as a {\em Bayes Decision} $BD(c,\pi)$ (the terminology is taken from \cite{brummer}). For any $\pi$, we denote the expected costs of the BD according to \eqref{expc} by $c(\pi)=c_\pi(BD(c,\pi))$.

We are now interested whether additional information will, in expectation, lower the costs of the BD. When additional information $e$ comes in, we will update $\pi$ to $\pi \mid e$, which we will write as $\pi(e)$. In this notation, $\pi(\emptyset)=\pi$, making clear that $\pi$ is the probability distribution  not conditioned on (any derivative of) the data $(e_x,e_y)$. We then take a Bayes Decision on the basis of $\pi(e)$. 
It may of course happen that this Bayes Decision happens to be, for some instances of $e$, more costly than the one based on $\pi$.
But since $\pi(e)$ is a better informed probability distribution than $\pi$, we would expect that on average, it is advantageous to take $e$ into account. We will show that this is correct in the sense that the following holds.
\begin{theorem}
	\label{bdimp}
Let $c=(c_{ij})$ be a cost function as described above for mutually exclusive and exhaustive hypotheses $H_1,\dots,H_n$ and actions $A_1,\dots,A_m$.  Let $\pi$ be the prior probability distribution on the $H_i$ and let $\pi(E)$ be the (random) posterior probability vector obtained from (a Bayesian update of the prior), with the random variable $E$ modeling the evidence. Then, we have 
\[
\BE[c(\pi(E))] \leq c(\pi),\]
where the expectation is over the evidence we obtain.
\end{theorem}
%where $E$ denotes the random variable that models the process of obtaining evidence so that the expectation is over the evidence that one obtains. 
We note that in the above theorem and throughout the paper, whenever we take an expectation, this is always over the capitalized random quantities in the expression.

Theorem \ref{bdimp} means that the average cost of a Bayes Decision will either remain the same, or decrease when we have updated the prior probability distribution $\pi$ to the posterior probability distribution $\pi(e)$. We prove this in the next sections. Note that Theorem \ref{bdimp} implies, for cost functions, all the inequalities (formulated there for scoring rules) in \cite{vergeer}, where the arguments were presented on a case-by-case basis, comparing different types of LRs (feature/score and common/specific source) to each other from a Bayes decision perspective using proper scoring rules. We come back to this in Section \ref{sec_proper}.

Since Theorem \ref{bdimp} is completely general, it also applies to $g(E)$, so that one can write

\[\BE [c(\pi(g(E)))] \leq c(\pi).\]

Also, Theorem \ref{bdimp} can be used multiple times. For example, when we consider two pieces of evidence $e_1$ and $e_2$, the costs of a BD will in expectation improve when we first take one of the $e_i$ into consideration, and then again when we also incorporate the other one. This means that Bayes Decisions based on the scores are an improvement over not incorporating anything, and also that given the score we again expect further improvement when we next consider the whole of the evidence. For any specific case, however, it may still be true that the Bayes Decisions become more costly when we do this, compared to when we stop at incorporating the score. All this is in sharp contrast to the conclusions drawn in \cite{neumann}, namely that score-based LRs should not be used.

We next give an elementary proof of Theorem \ref{bdimp}. We treat the case with two hypotheses and two actions separately: its proof is so particularly simple we do not want to withhold it, even if is also covered by the general case.

\subsection{Two hypotheses and two actions}
\label{sec_first}
We consider two mutually exclusive hypotheses $H_1$ and $H_2$ with $\BP(H_1)+\BP(H_2)=1$, and actions $A_1$ and $A_2$. We assume that  $c_{11}=c_{22}=0$, so that $A_i$ is the `correct' action if $H_i$ is true and can be executed without costs. Since Bayes decisions will remain the same when all costs are multiplied by the same factor, we can without loss of generality assume that $c_{21}=1$ and $c_{12}=\gamma >0$. The expected cost of $A_1$ is then equal to $\gamma \BP(H_2)$, and the expected cost of $A_2$ is equal to $\BP(H_1)$. Hence, if we were to decide at this point, one chooses $A_1$ if and only if
$$
r:=\frac{\BP(H_1)}{\BP(H_2)} \geq \gamma,
$$
that is, if the prior odds $r$ are at least $\gamma$. We are now in the same setting as Section 3 of \cite{vergeer} with the constant $\gamma$ playing the role of the threshold value $Th$.

Now consider that we are able to obtain evidence (data) to be denoted $e$. A realization of evidence $e$ leads to posterior probabilities  $\BP(H_1 \mid e)$ and $\BP(H_2 \mid e)$. We are interested in the expected costs of the BD based on the posterior probabilities. 

First, assume that $r < \gamma$ so that on the basis of the prior the BD is $A_2$, with expected costs $\BP(H_1)$. The new evidence $e$ will lead us to change the BD if the posterior odds become at least $\gamma$, i.e., when the likelihood ratio $\BP(e\mid H_1)/\BP(e \mid H_2)$ is at least $\gamma/r$. We write
\begin{equation}
\label{changep}
s:=\BP\left(LR(e) \geq \frac{\gamma}{r} \mid H_1\right),
\end{equation}
for this probability under $H_1$. 
Furthermore, it is known that (cf.\ the proof of Proposition 2.4.2 in \cite{meesterslootenbook})
\begin{equation}
\label{eigenschap}
\BP\left(LR(e) \geq t \mid H_2\right) = \BP(LR(e) \geq t \mid H_1) E(LR(e)^{-1} \mid LR(e) \geq t, H_1),
\end{equation}
from which we conclude, using \eqref{changep}, that
$$
\BP\left(LR(e) \geq \frac{\gamma}{r} \mid H_2\right) \leq s \cdot \frac{r}{\gamma}.
$$
If $H_1$ is true, the probability that the posterior odds exceed $\gamma$ is $s$. If that happens, we choose $A_1$ with no costs. If the posterior odds do not exceed $\gamma$, then we choose $A_2$ with costs 1, something which happens with probability $1-s$. Hence, under $H_1$ the expected costs are $1-s$.

If $H_2$ is true, we only make costs if the posterior odds exceed $\gamma$, and this happens with probability at most $s \cdot r/\gamma$. 
It follows that the expected posterior costs of the BD are at most
$$
\BP(H_1) (1-s) + \BP(H_2)\gamma \cdot s \cdot \frac{r}{\gamma} = \BP(H_1).
$$
Since the prior expected costs are $\BP(H_1)$, we see indeed that the expected costs of the BD are at most the same as based on the prior probabilities. A similar reasoning holds when $r \geq \gamma$. This proves Theorem \ref{bdimp} for this case. 

\subsection{The general case}
\label{sec_general}
We now return to the general case, dropping all assumptions on the costs $c_{ij} \in \BR$. As in the introduction, we allow for $n$ hypotheses $ H_1, \ldots, H_n$ and $m$ possible actions $A_1, \dots, A_m$, together with an initial probability distribution $\pi=(\BP(H_1), \ldots, \BP(H_n))$ representing our current conviction or knowledge. 

To prove Theorem \ref{bdimp} we start with a general observation.
Let, for $k=1, 2, \ldots$, $p_k = (p_{k,1}, \ldots, p_{k,n})$ be probability vectors, and let $t_1, t_2, \ldots $ be non-negative numbers such that $\sum_{k=1}^{\infty} t_k =1$. Then
$\sum_k^{\infty} t_k p_k$ is again a probability vector, and we have
\begin{equation*}
\label{hier}
	\begin{split}
		c \left( \sum_{k=1}^{\infty} t_k p_k \right) 
				&= \min_i \sum_{j=1}^n c_{ij} \sum_{k=1}^{\infty} t_k p_{k,j} \\		
				&= \min_i \sum_{k=1}^{\infty} \sum_{j=1}^n  c_{ij} t_k p_{k,j}  \\
				&= \min_i \sum_{k=1}^{\infty} t_k \sum_{j=1}^n  c_{ij} p_{k,j}  \\
				&\geq \sum_{k=1}^{\infty} t_k \min_i \sum_{j=1}^n  c_{ij} p_{k,j}  \\
				&= \sum_{k=1}^{\infty} t_k c(p_k).  \\
	\end{split}
\end{equation*}

Since
\[\BE [c(\pi(E))] = \sum_k \BP(e=e_k) c(\pi(e_k)),\]
taking $t_k = \BP(e=e_k)$ and $ p_k = \pi \mid e_k$ above leads to
\begin{equation}
%E_e (c(\pi \mid e)) \leq c \left( E_e(\pi \mid e) \right).
\BE [c(\pi(E))] \leq c \left( \BE[\pi(E)] \right).
\end{equation}
Alternatively, one can apply Jensen's inequality to arrive at the same conclusion.

It remains to show that $\BE[\pi(E)]=\pi$, i.e., that the expectation of the posterior probability distribution $\pi(e)$ is equal to the prior $\pi$. But this follows immediately from the fact that
$$
\BE(\BP(H_i \mid E)) = \sum_e \BP(H_i \mid e)\BP(e) = \BP(H_i),
$$
for all $i =1, \ldots, n$. This proves Theorem \ref{bdimp}.

This information-theoretical perspective is the key for understanding the distinction between the common source versus the specific source scenario on the one hand, and the difference between score-based and feature-based on the other. In both comparisons, the issue is that one of the alternatives is based on less information than the other. Theorem \ref{bdimp} expresses that on average, the better-informed situations leads to better decisions. That does not mean that the less informed LRs are `wrong' in whatever sense, only that they are based on less information. If this, however, is all the information that is available, then using this information is the right thing to do.

So far, these assessments have all been qualitative. Full features are better than scores, but to what extent? Or put differently, how bad is it to use less information?
Measuring performance is often a very contextual matter, but in general it is well known \cite{meesterslootenbook} that, when comparing hypotheses $H_1$ and $H_2$,
\begin{equation} \label{H0bound}\BP(LR \leq 1/t \mid H_1) \leq 1/t,\end{equation}
and
\begin{equation} \label{H1bound}\BP(LR \geq t \mid H_2) \leq 1/t.\end{equation}

Now, suppose that we consider two hypotheses $H_1$ and $H_2$ with some prior $\pi$ and data $e$, and that we have obtained $\pi(g(e))$ by computing $LR(g(e))$ for some $g(e)$.
Then, we can update $\pi(g(e))$ to $\pi(e)$ by calculating $LR(e \mid g(e))$, so that the two inequalities above apply to the computation of the LR for $e \mid g(e)$. Hence, if $H_1$ were true, then the probability that $LR(e \mid g(e))$ is more than a factor $t$ smaller than $LR(g(e))$, is bounded by $1/t$, and so is the probability that $LR(e \mid g(e))$ is more than a factor $t$ larger than $LR(g(e))$ if $H_2$ is true.

Hence, if for example $g(e)$ provides strong evidence for either hypothesis, then the probability that taking the full data $e$ into account will incorrectly (or rather, unfortunately) point strongly in the other direction is bounded as just described. We next illustrate these observations with an extensive example.

\subsection{Example: DNA kinship LRs}

In this example, the goal is to investigate whether two persons are siblings or unrelated. To that end, we consider the DNA profile of a person as a `trace measurement' of the DNA of their parents. Thus, the population of sources is the population of pairs of man and women (assumed to be always unrelated to each other), and sources give rise to measurements in the form of a DNA profile of a child of theirs. This model is well understood: DNA-profiles of parents are described by a probability distribution on the DNA profiles of persons in the general population, and children are obtained by Mendelian inheritance (perhaps enriched with a mutation model).

Now, suppose we have two individuals $X$ and $Y$ who are either siblings or unrelated. Write $H_0$ for the hypothesis that they are full siblings, and $H_1$ for the hypothesis that they are unrelated. Since we assume that these two hypotheses are exhaustive (i.e., have total probability equal to one), this example directly connects to Section \ref{sec_first}. Furthermore, we assume the following types of measurements: $e_x$ (the DNA profile of $X$), $e_y$ (the DNA profile of $Y$), and, say, $\theta_x$ (the DNA profiles of the parents of $X$).

If we only know $(e_x,e_y)$ we will evaluate these in our LR which is then equal to $\BP(e_x, e_y \mid H_0)/\BP(e_x, e_y \mid H_1)= \BP(e_y \mid H_0, e_x)/\BP(e_y \mid H_1, e_x)$. The computation for $H_0$ conceptually amounts to integration over all possible parents. Generally, the resulting LR will be larger when the profiles $e_x$ and $e_y$ share more alleles.
Now, consider that $\theta_x$ also becomes known, i.e., we get to know the DNA profiles of the parents of $X$. In that case, the profile $e_x$ becomes redundant: we can directly compare $e_y$ to its possible parents. The LR becomes $\BP(e_x, e_y, \theta_x \mid H_0)/\BP(e_x, e_y, \theta_x\mid H_1)  = \BP(e_y \mid H_0, \theta_x)/\BP(e_y \mid H_1, \theta_x)$.

Clearly, the LR based on $(e_x, e_y, \theta_x)$ is not the same as for $(e_x, e_y)$. In fact, it is easy to construct examples where the former is zero and the latter large, because there certainly exist DNA profiles $\theta_x, e_x, e_y$ such that $e_x, e_y$ are similar and lead to a large LR for being siblings, but also such that the parents of $X$ cannot be the parents of $Y$.

But this is absolutely no reason to dismiss, or in any way distrust, the LR based on the profiles of the alleged siblings only. When only $(e_x, e_y)$ are available, the evidence for being siblings versus unrelated is given by $LR(e_x, e_y)$. Note also that when this LR becomes larger, the posterior probability for being siblings increases as well, and if $X$ and $Y$ are indeed siblings then the scenario that we just sketched (additionally obtaining the profiles of the parents of $X$ and excluding the relationship with $Y$) is not possible.

Now, instead of processing the profiles $(e_x, e_y)$ we could consider a score $g(e_x, e_y)$, for example counting the number of alleles that the profiles have in common. Again, if only this score were known to us, the evidence is then $LR(g(e_x, e_y))$. There is nothing inherently wrong with this LR, but since the pair $(e_x, e_y)$ carries more information, we would rather assign a LR based on the full data, than only on the score.

To illustrate this, we have run a simulation experiment, generating all data just described: both for full siblings and unrelated individuals, we count the number $g(e_x,e_y)$ of shared alleles between $e_x$ and $e_y$, we compute the likelihood ratio $LR(e_x, e_y)$ comparing being siblings to being unrelated, and we compute the LR for $e_y$ to be a child of the parents of $X$, versus unrelated to both of them. We have done so based on Dutch allele frequencies (\cite{FLDOfreqs}) and using various DNA multiplexes: one with 10, and one with 15 autosomal loci. In all cases the profiles were simulated based on 15 loci and then calculations were done for these, as well for the subset of 10 loci. We simulated 500,000 cases. Note that these are, nowadays, rather small numbers of loci; we chose them purely for illustration purposes.

\begin{figure}
	\centering
	\begin{subfigure}{0.48\textwidth}
		\centering
		\includegraphics[width=0.95\textwidth]{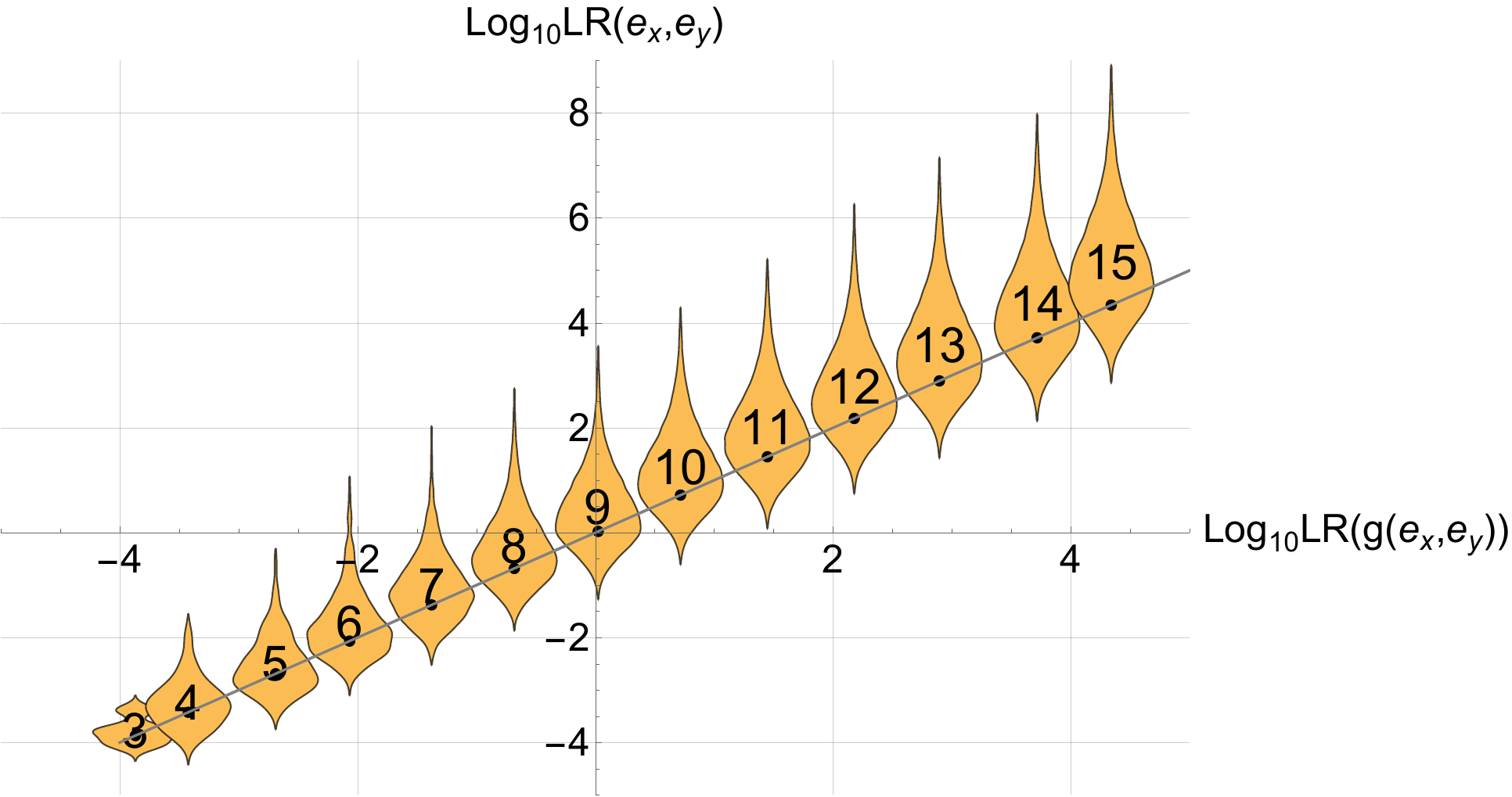} % first figure itself
		\caption{true relationship: siblings}
        \label{CountVersusSI10sibs}
	\end{subfigure}\hfill
	\begin{subfigure}{0.48\textwidth}
		\centering
		\includegraphics[width=0.95\textwidth]{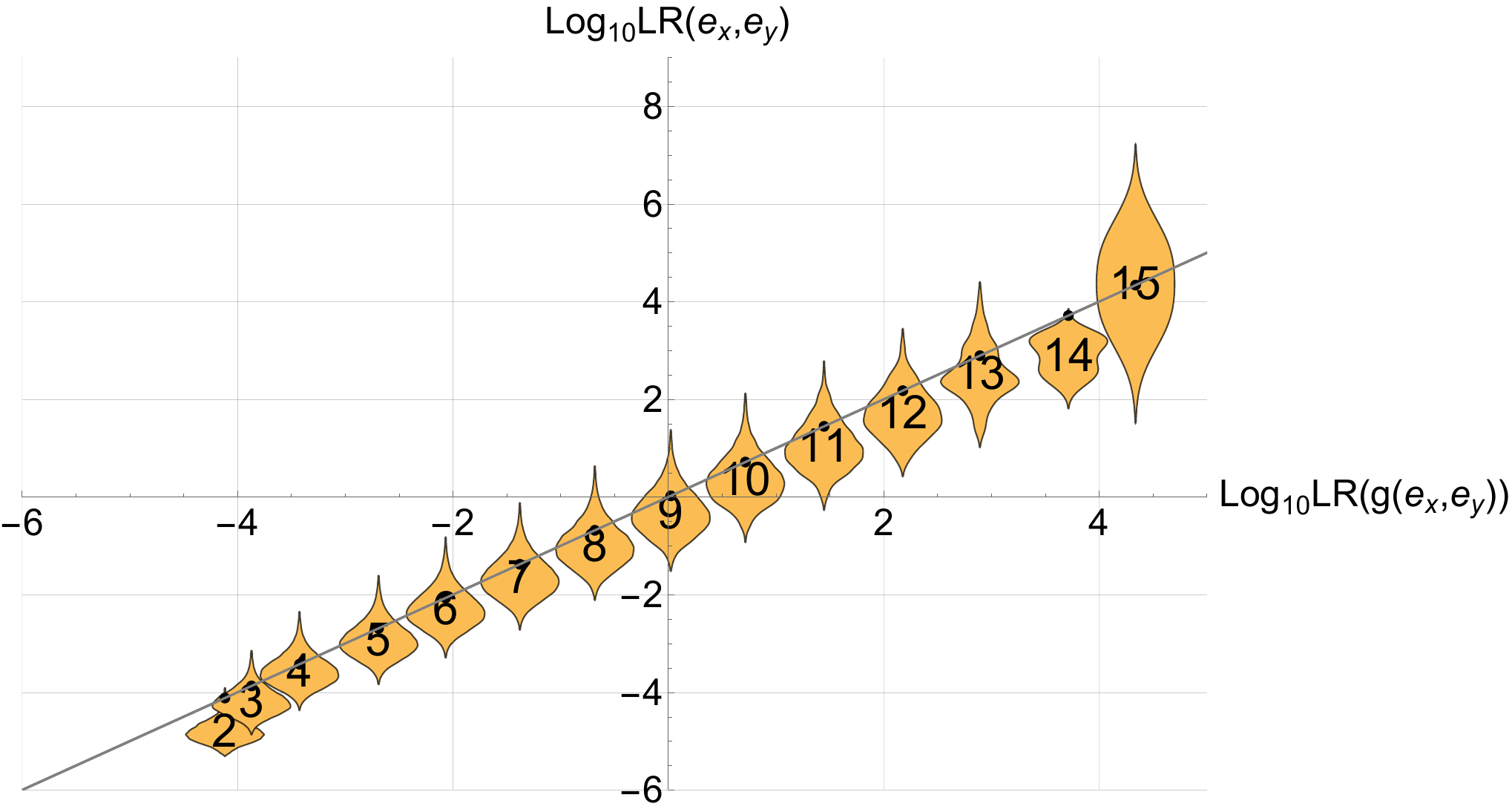} % second figure itself
		\caption{true relationship: unrelated}
        \label{CountVersusSI10unrel}
	\end{subfigure}
    \caption{$\log_{10}(LR(e_x,e_y))$ (based on DNA profiles) versus $\log_{10}(LR(g(e_x,e_y)))$ (based on number of shared alleles). Each violin plot represents the $LR(e_x,e_y)$ for profiles whose number of shared alleles is displayed in the plot.}
    \label{CountVersusSI}
\end{figure}

First of all, suppose that we work on 10 loci. We obtain likelihood ratios $LR(g(e_x,e_y))$ for being siblings, versus unrelated, based on the empirical number $g(e_x,e_y)$ of alleles shared on 10 loci. We can compare these to the LRs obtained when $LR(e_x,e_y)$ is computed, and the result is displayed in Figures \ref{CountVersusSI}.

In these figures, we have placed violin plots representing the distribution of $LR(e_x,e_y)$ placed at the values $LR(g(e_x,e_y))$  corresponding to a number $g(e_x,e_y)$ of shared alleles, and we see for example that neutral evidence is obtained for 9 (out of possibly 20) shared alleles. We see in these figures that, when the profiles are evaluated, the $LR(e_x,e_y)$ are of course different from $LR(g(e_x,e_y)) $,  and have a tendency to more strongly support the correct hypothesis. It may, of course, happen that the LR based on $(e_x,e_y)$ supports the other hypothesis compared to the LR based on $g(e_x,e_y)$.  But we also see that this becomes less likely for larger $|Log_{10}LR(g(e_x,e_y)|$, as predicted. Finally we note that each violin plot has the same width irrespective of the number of instances of $g(e_x,e_y)$.

\begin{figure}
	\centering
	\begin{subfigure}{0.45\textwidth}
		\centering
		\includegraphics[width=0.9\textwidth]{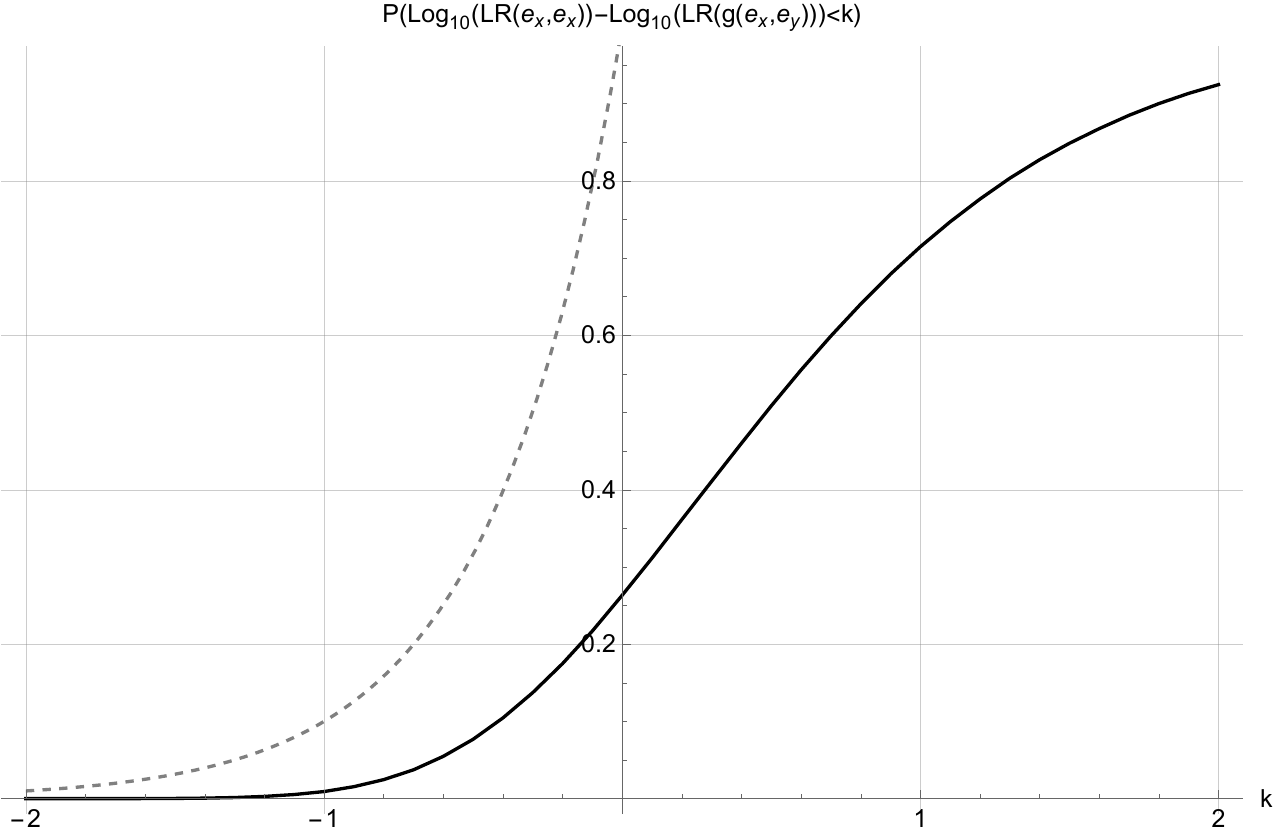} % first figure itself
		\caption{Actual siblings}\label{CountVersusSI10sibsCDF}
	\end{subfigure}\hfill
	\begin{subfigure}{0.45\textwidth}
		\centering
		\includegraphics[width=0.9\textwidth]{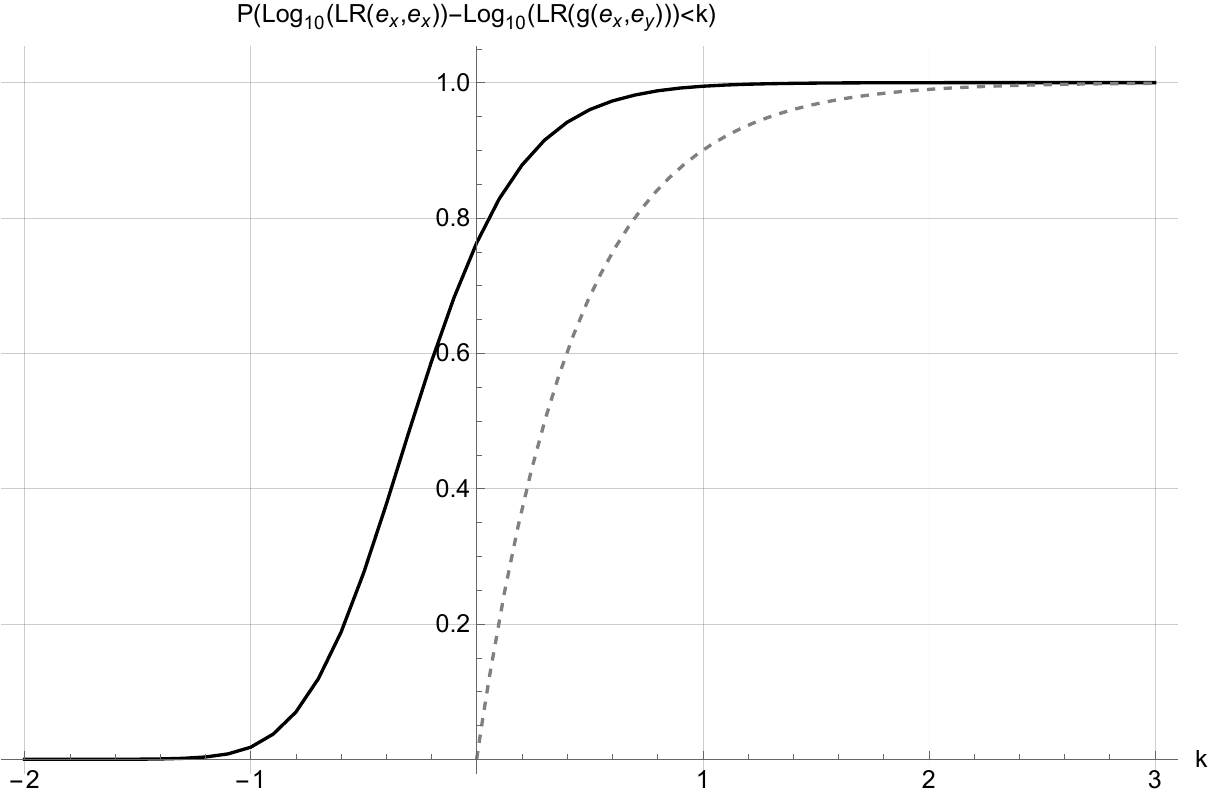} % second figure itself
		\caption{Unrelated pairs}\label{CountVersusSI10unrelCDF}
	\end{subfigure}
    \caption{Cumulative distribution function of $\log_{10}(LR(e_x,e_y))-\log_{10}(LR(g(e_x,e_y)))$ (black), and theoretical bound (dashed).}
    \label{CountVersusSI10CDF}
\end{figure}

In Figure \ref{CountVersusSI10CDF} we plot the distribution of the difference between the two LRs on a logarithmic scale. Indeed, comparing $LR(e_x,e_y)$ with $LR(g(e_x,e_y))$ by considering their quotient, we see that  $LR(e_x,e_y)/LR(g(e_x,e_y))$ respects the bounds \eqref{H0bound} and \eqref{H1bound}.  Clearly, both for siblings and unrelated individuals the evidence tends more strongly towards the hypothesis that is actually true. There are, also, cases where this does not happen, which is inevitable and to be expected. The frequency with which this occurs, however, is bounded by \eqref{H0bound} and \eqref{H1bound}.

\begin{figure}
	\centering
	\begin{subfigure}{0.45\textwidth}
		\centering
		\includegraphics[width=0.9\textwidth]{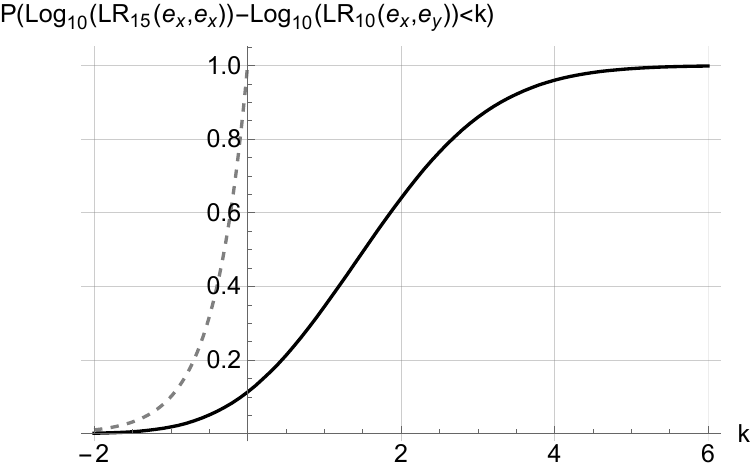} % first figure itself
		\caption{Actual siblings}\label{SI15vs10sibs}
	\end{subfigure}\hfill
	\begin{subfigure}{0.45\textwidth}
		\centering
		\includegraphics[width=0.9\textwidth]{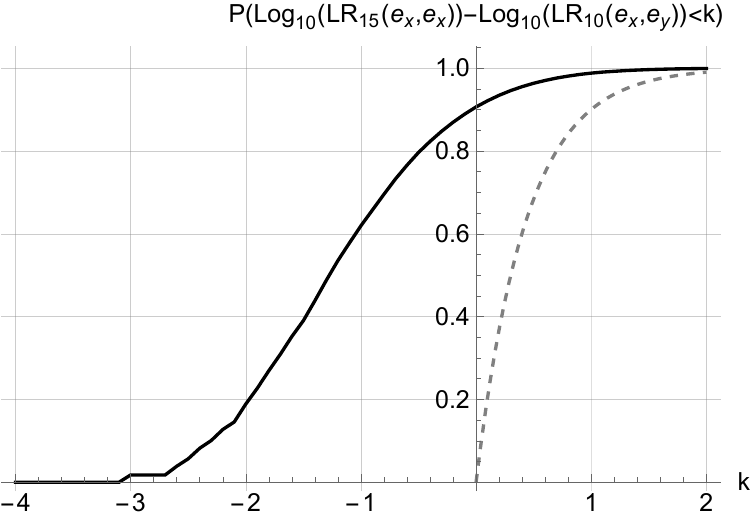} % second figure itself
		\caption{Unrelated pairs}\label{SI15vs10unrel}
	\end{subfigure}
    \caption{Difference in $\log_{10}(LR(e_x,e_y))$ on 15 versus 10 loci.}
\end{figure}

Thus, we see that when we process the profiles $(e_x, e_y)$ instead of the number of shared alleles $g(e_x,e_y)$, the LRs change. That is no reason to say that, in hindsight, the LR based on $g(e_x,e_y)$ is incorrect; it is only a sub-optimal way to treat the data. But so is $(e_x, e_y)$ when more loci become available: comparing the LRs for $(e_x, e_y)$ on 15 with those restricted on 10 loci gives similar changes, as shown in Figures \ref{SI15vs10sibs} and \ref{SI15vs10unrel}. Note that, since for these loci the genetic data of unrelated persons are independent on the loci considered, and so is the inheritance of alleles towards offspring, these graphs also represent the likelihood ratio distribution on the five additional loci.

\begin{figure}
	\centering
	\begin{subfigure}{0.45\textwidth}
		\centering
		\includegraphics[width=0.9\textwidth]{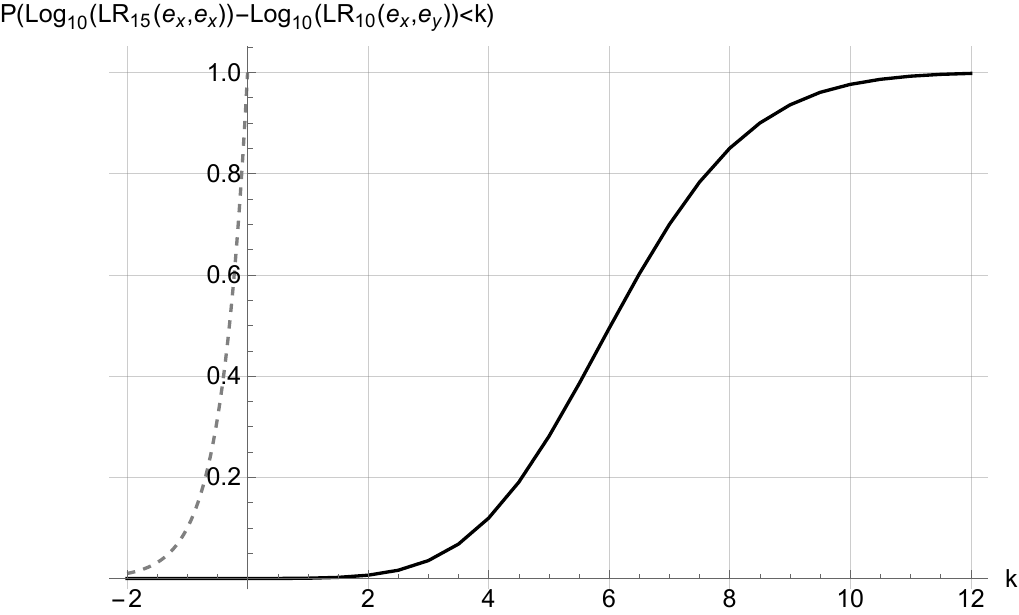} % first figure itself
		\caption{Actual siblings. }
        \label{LRparVersusSI10sibs}
	\end{subfigure}\hfill
	\begin{subfigure}{0.45\textwidth}
		\centering
		\includegraphics[width=0.9\textwidth]{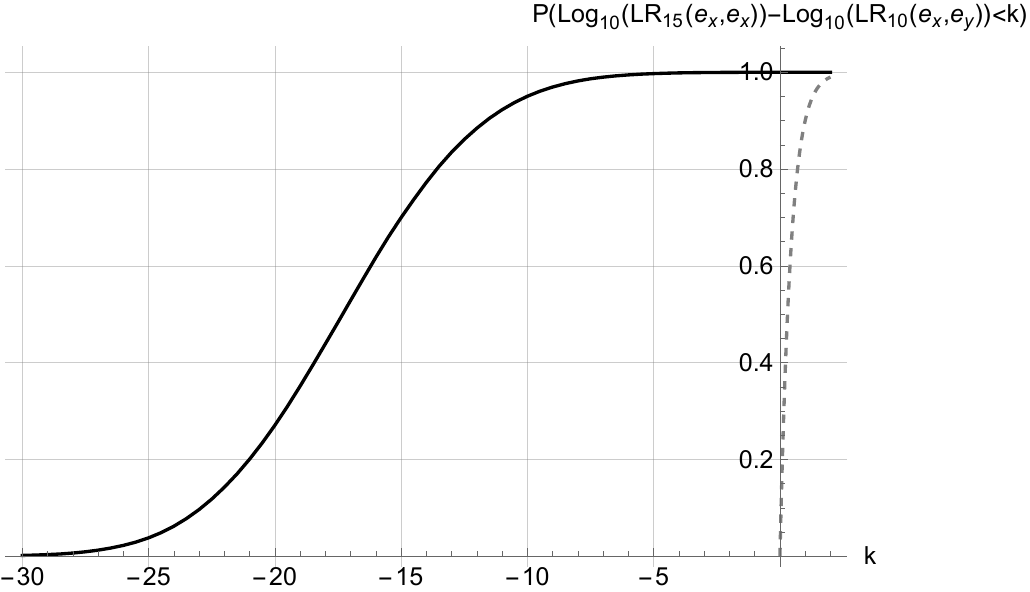} % second figure itself
		\caption{Unrelated pairs}
        \label{LRparVersusSI10unrel}
	\end{subfigure}
    \caption{Difference in $\log_{10}(LR)$ based on profiles, with or without parents of $X$.}
    \label{LRparVersusSI10}
\end{figure}

The most dramatic changes in LR occur, of course, if the parents of $X$ become known and the question of whether $X$ and $Y$ are siblings reduces to whether $Y$ is a child of the parents of $X$. In that case, going back to 10 loci, we get the changes in LR displayed in Figure \ref{LRparVersusSI10}.

We see that the extra information contained in the profiles of the source (i.e., the parents of $X$) gives much stronger likelihood ratios. Still, this is no reason to disqualify LRs based on only the profiles of $X$ and $Y$ alone, if this is all that is available. 

In conclusion, we see in this example how the different ``LR systems'' can all be defined in a context where they are not traditionally thought of as such; indeed, all these LRs simply differ in which information is available, and whether the full information is processed or only a function of it. They all compare the hypothesis ``the source of $Y$ is $S_x$, the source of $X$" to the alternative ``the source of $Y$ is a random unknown source", either on the basis of $g(e_x,e_y)$, $(e_x,e_y)$ or $(e_x, e_y, \theta_{x})$. Of course, when such a calculation is done, it is of utmost importance that the forensic report contains a full description of the considered hypotheses and the processed data. A LR based on the profiles of the alleged siblings alone, ignoring the profiles of the parents, is not wrong: it expresses the information obtained from what is taken into account. At the same time, of course, if more data is available for which a statistical model is available, any analysis that ignores these data is sub-optimal. Although we assumed only two options throughout this example, the general conclusions also hold in case of more or two non-exhaustive hypotheses. In those cases, the LRs calculated above only address the strength of evidence relative to two hypotheses at hand. The point remains that more information will, on average, never yield worse decisions.

\section{An analysis of some arguments against score-based methods in a toy example}
\label{sec_neumann}

We have discussed common and specific source LRs, as well as the distinction between a feature-based and a score-based LR. We saw that each differs in what information it processes, and is valid in the sense that if the model is accurate, then the expected costs of Bayes Decisions can only decrease when we incorporate more information. Thus, depending on amount of information available, there is no reason to reject some of them on principle grounds. 

In \cite{neumann}, the opposite conclusion was drawn. Therefore, this section is devoted to a critical analysis of their arguments.

\subsection{A toy model}\label{toy}
The goal in \cite{neumann} is to ``separate the wheat from the chaff'' when calculating likelihood ratios in various ways. To that end, they set up a toy model that allows for relatively simple computations. The setup is such that the computations can be done in several ways: score/feature-based as well as common/specific source.
They consider two physical sources $S_x$ and $S_y$. The goal is to evaluate forensic evidence $e=(e_x,e_y)$ coming from these sources with respect to the following two hypotheses:
\begin{itemize}
	\item $H_0$: $S_x=S_y$,
	\item $H_1:$ $S_x \neq S_y$.
\end{itemize}
The toy model in \cite{neumann} is as follows (cf.\ their section 2.3).  Any source is represented by a real number that is the expected value of some hypothetical measurement of (a property of) that source. On the population of sources, this expected value is normally distributed with mean $\mu$ and variance $\sigma_D^2$. When we measure on a source, we obtain a random perturbation of its expected value. If that expected value is $\mu_d$ (this corresponds to $\theta_S$ for source $S$ in our previously introduced notation, but here we keep the notation of \cite{neumann}), the outcomes for measurements on this source are modeled by a normal distribution with mean $\mu_d$ and variance $\sigma_s^2$ (when the measurement is modeled as a reference measurement) or variance $\sigma_u^2$ (when the measurement is modeled as a trace measurement). The distinction between $\sigma_u$ and $\sigma_s$ allows to take into account that measurements on a reference sample can be performed under more ideal conditions and then tend to give results more closely to the actual value $\mu_S$ of the source $S$ they come from.

Next, still following \cite[2.3]{neumann}, a common source scenario consists of two measurements $(e_{u_1}, e_{u_2})$. Both observations involve $\sigma_u^2$ now, since two traces are modeled. The authors write $H_{0,cs}$ and $H_{1,cs}$ for the hypotheses that the measurements are on the same, or different sources, and write the resulting likelihood ratio as 
$$
LR_{CS}=\frac{\BP(e_{u_1}, e_{u_2} \mid H_{0,cs})}{\BP(e_{u_1}, e_{u_2} \mid H_{1,cs})}=\frac{f(e_{u_1}, e_{u_2} \mid H_{0,cs})}{f(e_{u_1} \mid H_{1,cs})f(e_{u_2} \mid H_{1,cs})},
$$
the latter expression using that different sources are modeled as having independent parameters. The function $f$ represents the normal density with appropriate parameters.
	
This is contrasted to a specific source scenario, in which we have measurements $e_{u}$ and $e_s$. Here, $e_{u}$ is a trace measurement as before, and $e_s$ is a reference type measurement on a known source $S$ , i.e., with known parameter $\mu_d$. The authors write $H_{0,ss}$ for the hypothesis that $e_u$ has source $S$, and $H_{1,ss}$ for the hypothesis that the source of $e_u$ is unknown.  Next, they claim that
$$
LR_{SS}=\frac{\BP(e_u, e_s \mid H_{0,ss})}{\BP(e_u, e_s \mid H_{1,ss})}=\frac{\BP(e_u \mid H_{0,ss})}{\BP(e_u \mid H_{1,ss})}. 
$$
However, this is only correct when $e_u$ is independent of $e_s$, given $H_{\cdot,ss}$. For this to hold, we must also condition on the source value $\mu_d$, otherwise $e_s$ is informative for the source. If we do condition on $\mu_d$, the observation $e_s$ is redundant, representing a superfluous measurement on a source whose characteristic $\mu_d$ is already known, leading to the above equality. In other words, if $\mu_d$, the parameter of source $S$, were not conditioned on, the LR would become 
$$\frac{\BP(e_u \mid H_{0,ss}, e_s)}{\BP(e_u \mid H_{1,ss}, e_s)},
$$
and $e_s$ cannot be omitted from the conditioning. Hence, \cite{neumann} tacitly assume knowledge of $\mu_d$ in their definition of the hypotheses $H_{\cdot,ss}$. Apparently, while the specific source and common source hypotheses are very similar in notation, those for the specific source include the value $\mu_d$  of the source of $e_s$ whereas the hypotheses for the common source LR do not contain such a value, since in this setup there is no (measured) candidate source for any of the measurements.  We will see below that this clarification of the notation is key in understanding the logical errors in the obtained conclusions of \cite{neumann}.

In their section 2.4, the authors set up a comparison between LRs obtained with the common source framework and with the specific source framework, seemingly setting up a comparison where $LR_{CS}$ and $LR_{SS}$ as above are computed on the same data: first, data $(\mu_d, e_u, e_s)$ are generated. Then, the data $(e_u, e_s)$ are used for a common source evaluation. The authors state that ``To calculate the common source likelihood ratio using the data generated under the specific source model, we set $e_{u_1}=e_u$; $e_{u_2}= e_s; \sigma_{u_1}^2=\sigma_u^2$ and $\sigma_{u_2}^2=\sigma_s^2$". Thus, the simulation results seem to be presented as though the same data are evaluated by different models. This is, even, explicitly stated in section 2.2. preceding the toy model where we read, about the difference between common source and specific source: ``Each scenario results in {\it different likelihood functions for the same information}, and in different interpretations of the results of forensic examination" (our italics). 

This however is inaccurate, since for $LR_{SS}$ not only $(e_u, e_s)$ are evaluated, but the triples $(e_u, e_s, \mu_d)$ in which $\mu_d$ is the parameter of the source that $e_s$ comes from, and according to $H_{0,ss}$ also $e_u$. Again, $e_s$ is redundant in this triple. Thus, the actual comparison is between $(e_u, \mu_d)$ and $(e_u, e_s)$, both for $H_0$ versus $H_1$. In order to compare the resulting LRs, one should randomly sample $(e_u, e_s, \mu_d)$ under the two hypotheses, and then compare the LR processing all relevant information $(e_u, \mu_d)$ with the situation in which one processes $(e_u, e_s)$ only. Instead, the authors choose three different cases: they fix $\mu=10, \sigma_D^2=10, \sigma_u^2=2$ and then consider the cases corresponding to $\mu_d=9, \sigma_s^2=1$, or $\mu_d=0, \sigma_s^2=1$, and finally $\mu_d=9, \sigma_s^2=10^{-5}$. In the last case, almost no differences between $LR_{SS}$ and $LR_{CS}$ are obtained, which is quite reasonable since the small $\sigma_s^2$ means that $e_s$ all but reveals $\mu_d$. The differences for the other cases lead the authors to dismiss the common source LRs. As an aside, note that the second case is quite irrelevant, since it represents a source with parameter about three standard deviations away from the mean; such sources will only very rarely be observed. 

However, these conclusions are not justified at all. As we explained, the authors are not, as they claim, evaluating the same data under different sets of hypotheses. They do precisely the opposite: they compare different data for the same hypotheses, where the question is whether two measurements are measurements on the same source or not. There is only one model, namely the toy model described above, that describes all possible data. In other words, this is not a matter of ``models'', ``systems'' or ``frameworks'', but rather a matter of available data. In the specific source framework, more information is available than in the common source framework, but the statistical model and the hypotheses are the same. Furthermore, the comparison is not a random selection of evidence for evaluating the same versus different sources, and the simulation data are therefore not suitable for general conclusions.

The simulations in \cite{neumann} do not allow to infer probabilistic assessments of the distribution of specific-source LRs given a certain specific-source LRs, due to the nature of the simulation.
Their simulations compare $LR(\mu_d, e_u, e_s)$ with $LR(e_u, e_s)$ for a fixed choice of $\mu_d$. However, in order to judge how much these LRs differ in practice, we would need a probability distribution describing what $LR(\mu_d, e_u, e_s)$ may result, given $(e_u, e_s)$ or given $LR(e_u,e_s)$. That is, we would need to update the probability distribution for the source of $e_u$ with the two measurements via $LR(e_u, e_s)$, in order to then sample $\mu_d$ from the updated probability distribution and obtain the distribution for $LR(\mu_d, e_u, e_s)$.  

\begin{figure}
	\centering
	\begin{subfigure}{0.5\textwidth}
		\centering
		\includegraphics[width=\textwidth]{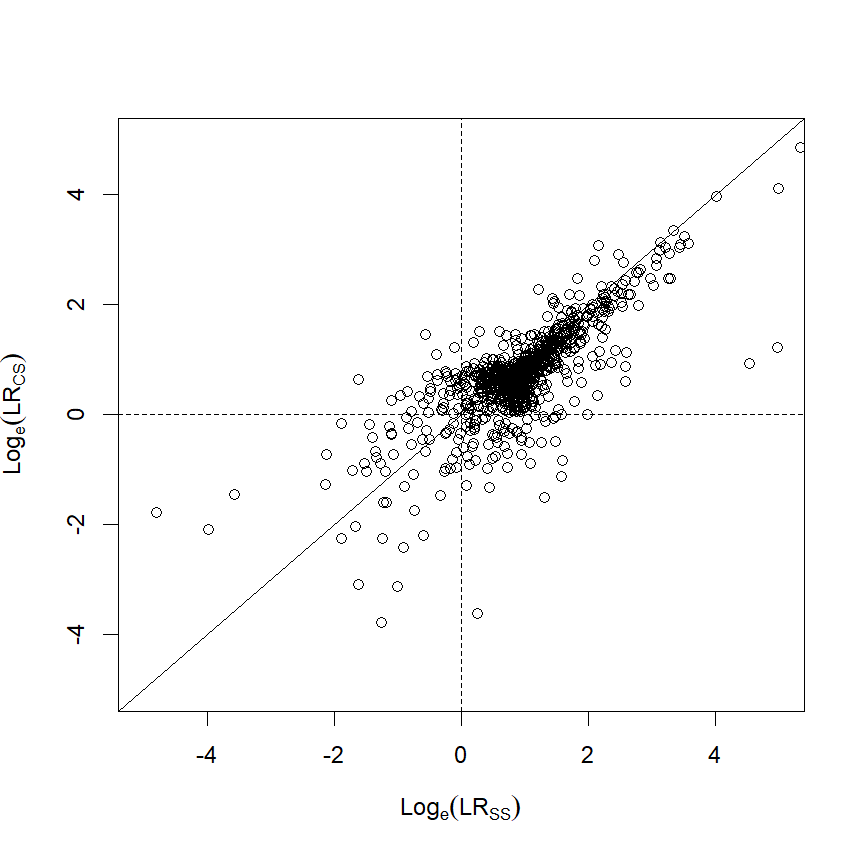} % first figure itself
		\caption{$H_0$ true (same source)}
        \label{plaatje1}
	\end{subfigure}\hfill
	\begin{subfigure}{0.5\textwidth}
		\centering
		\includegraphics[width=\textwidth]{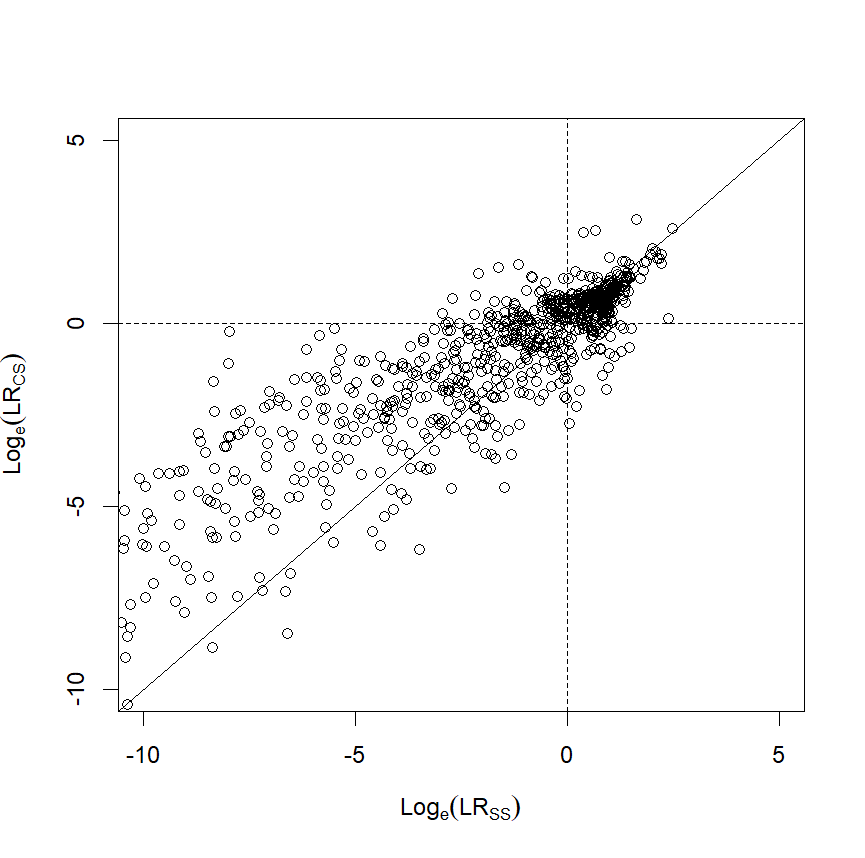} % second figure itself
		\caption{$H_1$ true (different sources)}
        \label{plaatje2}
	\end{subfigure}
    \caption{Comparison of common source likelihood ratios $LR_{CS}$ and specific source likelihood ratios $LR_{SS}$, obtained with the toy model.}
    \label{dinges}
\end{figure}

In contrast, here the value of $\mu_d$ is fixed so that no impression can be obtained from the distribution of specific-source LRs that one would obtain given an outcome of the common-source LR. The distributions in the simulations are, therefore, of limited value beyond studying what may happen for very specific parameter choices. For completeness, we carry out the required comparison, even if this is strictly speaking not needed since we already know from the general principles above that the data $(e_u, \mu_d)$ will on average be more informative than $(e_u, e_s)$, and that always only processing $(e_u, e_s)$ is also a perfectly reasonable procedure, if the full evaluation is impossible, for whatever reason. In \Cref{dinges}, simulation results are plotted with $\mu=10, \sigma_D^2=10, \sigma_u^2=2,\sigma_s^2=1,$ all values in accordance with the first two settings of \cite{neumann}. However, we let the source value $\mu_d$ vary in accordance with the population. We believe that this setting should have been used instead of picking specific values for $\mu_d$. As expected, the specific source information generally leads to better results, but the effect is on average not as dramatic as Figure 1 (b) of \cite{neumann} suggests.

Similar flaws occur in the reasoning leading to Figures 3 and 4. %The authors continue with an example of suspect-centered score-based likelihood ratios. Here, for the specific source scenario, they still evaluate $\mu_d, e_s$ and $e_u$ as before, and compare this to the  evaluation of $(e_{u}-e_{s})^2$ and $e_s$. As before, the results are presented as a comparison of the same data in two different frameworks, and as representative, whereas the opposite is true: different data are processed for the same hypotheses, and the simulations are by no means representative. The authors, again, conclude against the use of common source LRs: ``common source score-based models may be convenient to implement, but are not relevant to most examinations of forensic interest, and do not equal to the weight of forensic evidence."
Again, whereas their simulations results are accurate given their model and parameter choice, the conclusions are not justified.

\subsection{Lack of coherence?}
Finally, we comment on another claim in \cite{neumann}. In their Section 3.3, they discuss specific source score based likelihood ratios. Such a LR evaluates as evidence the score $\delta(e_u,e_s)$ and $e_s$ where $e_s$ represents, as before, a control measurement on the putative source. It is not clear from the notation whether the `true' value $\mu_d$ of this source is considered to be known as part of the hypothesis. This seems likely in view of it being so for $H_{\cdot,ss}$ discussed before. However, if $\mu_d$ is known it is clearly suboptimal to take the score $\delta(e_u,e_s)$  instead of $\delta(e_u,\mu_d)$, which might suggest that $\mu_d$ is not conditioned on here. %When the score $\delta(e_u,e_s)$ is evaluated as evidence, one considers the measurement $e_u$ w.r.t. $e_s$, and one cannot omit $e_s$ as superfluous measurement.
For the analysis below, whether the true value of the source is conditioned on is not important.
Now, they claim that this setup is ``plagued by a fundamental lack of coherence''. With that, they mean the following. Suppose that, instead of one known source, we now have two known sources $A$ and $B$ with respective parameters $\mu_a$ and $\mu_b$, and we still have a measurement $e_u$ on a trace.
Now two hypotheses $H_A$ and $H_B$ are considered with respect to $e_u$: it either is a measurement of source $A$ or of source $B$, corresponding to the two possibilities
\begin{itemize}
    \item $H_A$: $e_u \sim \CN(\mu_a,\sigma_a^2)$,
    \item $H_B:$ $e_u \sim \CN(\mu_b,\sigma_b^2)$.
\end{itemize}
In addition we have measurements $e_a$ and $e_b$ obtained from $A$ and $B$. From these, one computes the scores $\delta(e_u, e_a)$ and $\delta(e_u, e_b)$.
%The authors now consider two pieces of evidence. In the first case, they assume the value $\mu_A$ is known. Then the evidence is $\mu_A$, together with the score $\delta(e_u,\mu_a)$. In the second case, we know the value $\mu_b,$ yielding the score $\delta(e_u,\mu_b)$. 
Now \cite[(10)]{neumann} notes the fact that in general,
\begin{equation}
	\begin{split}
		\frac{\BP(\delta(e_u,e_a) \mid e_a, H_A)}{\BP(\delta(e_u,e_a) \mid  e_a, H_B)} \neq \frac{\BP(\delta(e_u,e_b) \mid  e_b, H_A)}{\BP(\delta(e_u,e_b) \mid  e_b, H_B)}.
	\end{split}
\end{equation}
They see this as showing that the specific source score based likelihood ratio application ``is not coherent in general since it potentially does not support the same proposition depending on which one is considered first".  The above inequality however only expresses that using different data may of course lead to a different LR. Here, either the pair $(\delta(e_u,e_a),e_a)$ is processed, or the pair $(\delta(e_u,e_b), e_b)$, but never the full data. Suppose that there is no measurement error, so that $\mu_a=e_a$ and $\mu_b=e_b$ are known. Then the above inequality simply says that if we evaluate only the distance from $e_u$ to $\mu_a$, we get another likelihood ratio than when we only evaluate the distance from $e_u$ to $\mu_b$, which is perfectly logical. 
If all data would be taken into account, we would arrive at % the order in which this is done does not matter. Indeed,
\begin{equation}
	\begin{split}
		\frac{\BP(\delta(e_u,e_a),\delta(e_u,e_b), e_a,e_b\mid H_A)}{\BP(\delta(e_u,e_a),\delta(e_u,e_b), e_a,e_b\mid H_B)} = \frac{\BP(\delta(e_u,e_a),\delta(e_u,e_b) \mid e_a,e_b,H_A)}{\BP(\delta(e_u,e_a),\delta(e_u,e_b)\mid e_a,e_b, H_B)} .
	\end{split}
\end{equation}
%In this expression, the information about $A$ and $B$ is jointly incorporated, so that the order in which it presented does not matter, since it trivially holds that $A \cup B = B \cup A$ for any two sets $A,B$. As before, whatever data one takes into account, one gets a likelihood ratio that expresses the evidential value of those data. %On a side note, even if the alleged difference caused by the ordering of $H_A$ and $H_B$ would have existed, it would have had nothing to do with the score-based nature of SSS. Regardless of scores, using different data will generally lead to a different outcome.

\section{Discussion and conclusions}

The central question of this paper is whether or not additional evidence is somehow beneficial. In Sections \ref{sec_first} and \ref{sec_general} we showed that this is the case for any cost function, at least in expectation. Even if the cost function's `goal' is to lie about the truth, then one can tell lies better with more information. 

All case-by-case inequalities for cost functions derived in \cite{vergeer} are contained (but noting that he only treated the $n=2$ case), for cost functions $c_{i,j}$, by Theorem \ref{bdimp}. For instance, consider the situation in which one compares features with scores. %As we explained in the introduction, it is not uncommon that data is summarized into scores, a process by which inevitably information is lost.
Vergeer showed for this case (and, separately, for others) that on average, using scores leads to higher costs than using all data.  These conclusions all directly follow from our results. Indeed, one can interpret the `prior'  $\pi$ as the information one has after evaluating the score only. This $\pi$ may be the Bayesian update of an earlier prior upon seeing the score. To take the full features into account in addition to the already known score can in our set-up be interpreted as gathering extra data, and this leads therefore to lower expected costs. The same reasoning goes through for all other examples discussed in \cite{vergeer}. 

This principle also shows that calculating a score-based LR, while less informed than using the full data, is on average still better than using only the prior (3.1.3 in \cite{vergeer}). This directly contradicts the statement made in \cite{neumann} that within the Bayesian paradigm, ``one cannot use score based likelihood ratios''. We have analysed the arguments from \cite{neumann} in some detail, as we find it important to explain how the authors of that paper arrived at their conclusion, and why this argument is incorrect.

We would like to add that `all data' is in a sense deceiving, since more knowledge may be possible in the future. In the past, DNA typing was carried out with fewer loci than nowadays. For example, a widespread multiplex was the SGMPlus kit which contained 10 autosomal loci \cite{meesterslootenbook}. More recent kits contain  (about 10 to 15) additional loci, so that calculating a likelihood ratio based on SGMPlus has become a score based likelihood ratio, as it comprises a data reduction relative to the larger sets of loci that are nowadays inspected. This example reveals that what we regard as a full feature method today, may be very well be a partial features method in the future, providing another argument against banning score-based methods.

In fact, every method (pathological examples aside) to calculate a likelihood ratio may in fact be regarded as score-based. One method may of course be more score-based than others, but all continuous measurement methods will usually use some form of simplification (e.g., rounding, binning, cleaning data) and thus lose information that could theoretically be incorporated.

The current paper, like \cite{vergeer} and \cite{neumann}, has mostly been concerned with fundamental properties of score-based methods. We hope that this work contributes in clearing score-based likelihood ratios of any fundamental suspicion, and that future research can focus on other important aspects such as calibration \cite{ferrer, hannig, ramos2021,rodriguez, ypma}. 

As a consequence, the question is not whether we should use score based methods, but how they can be used best. What are the conditions needed to justify data reduction, and how good are the decisions made? The answer to these questions depends entirely on the context and the objectives, see also \cite{garton}. Some score-based methods that use substantial reduction (quite naturally) perform badly \cite{enzinger}, and may thus not be generally recommended for use in legal contexts.

%We also take a position in the alleged distinction between common source and specific source methods \cite{saunders}. As we pointed out, it is a matter of available information which scenario applies, while the underlying statistical model remains the same. As with scores versus features, we see no dichotomy between these alleged different ``LR-frameworks''.

Next, we come back to the perceived difference between the common source LRs versus the specific source LRs. Suppose we have a trace measurement $e_u$ and a source measurement $e_s$ of source $S$, and that the question is whether $e_u$ also is from $S$ or not. Following \cite{neumann} we only speak of a specific source LR when the parameters $\theta_S$ of the source are known, which can be realized for discrete parameters, but will never really be the case for continuously distributed parameters. Hence, the LR is either calculated based on $(e_u, e_s, \theta_s)$ (and called specific source LR) in which case $
\theta_S=e_s$ or $e_s$ is redundant, or on $(e_u, e_s)$ (and called common source LR). We see no principle difference between these two approaches here: in both cases, the resulting LR involves (for $H_1$) to first probabilistically infer the source parameters (if $\theta_S$ is known, this is a trivial step) and then integrating the probability to see $e_u$ over the resulting probability distribution for $\theta_S$. In both cases, the hypotheses are the same, namely whether the source $S$ is the source of $e_u$ or not. We reiterate that the terminology (specific source or common source) only indicates the type of information that is processed, it does not indicate different statistical frameworks.

Similarly, if $e_{u_1}$ and $e_{u_2}$ are measurements on two traces, we may set up the `common source' hypotheses that they are from the same, versus different, sources. However, we may equally well view this as `specific source' LRs for the hypothesis that $e_{u_2}$ comes from $S_1$, defined as the source of $e_{u_1}$, or not. If additional data $e_{s_1}$ on this source are known, then we would ideally evaluate $(e_{u_1},e_{u_2},e_{s_1})$. However, when we process only $(e_{u_1}, e_{u_2})$, the LR that we obtain is still a LR for the hypothesis that $S_1$, the source of $e_{u_1}$ is also the source of $e_{u_2}$, albeit one that has ignored relevant information. 
%Analogously, consider that we have trace measurements $(e_{u_1}, e_{u_2})$, and source measurement $e_s$ on source $S$, known to be the source of $e_{u_1}$. Further we assume that the question at hand is whether trace $e_{u_2}$ also originates from $S$. If only a common source model for traces is available, we process only $(e_{u_1}, e_{u_2})$ and ignore $e_s$, which typically will lead to less discriminating LRs. 
In our DNA example, this would amount to testing persons $A$ and $B$ for being siblings, ignoring the DNA profiles of the parents of $A$, hence also considering the possibility that $A$ and $B$ might be siblings with other parents. Nevertheless, we do investigate whether the `source' of $A$ (their parents) are also the source of $B$.
% It still is the case though that the LR we get is pertaining to the question whether $e_{u_2}$ comes from $S$ or not, since $S$ is by definition the source of $e_{u_1}$.

In general, we conclude that the common source LR operates on a smaller set of evidence than the specific source LR, but addresses the same question. %Since the common source LR uses less data, more uncertainty remains on the common source. 
Common source LRs quantify  evidence for the hypotheses that the source of the first trace is (versus is not) also the source of the second trace, just a specific source LRs do, but based on less information. That does not make them wrong or redundant; while suboptimal, it is preferable to use some evidence than no evidence. 

If large posterior odds on the formulated hypotheses are obtained based on a subset of the evidence, say with a score-based LR and/or common source LR, then this means (if one of the hypotheses must be true) that $H_1$ is likely correct. If that is so, the evidence not analyzed thus far will yield a LR, when evaluated, that in expectation supports $H_1$ also. In any case, the probability that new evidence, when evaluated, will yield a LR of at least $t$ in the direction of the hypothesis that is false, is bounded by $1/t$.  It is generally of course not impossible that the evidence, unevaluated so far, can cause a large change in evidential value. The probability of that happening, in view of the bounds for misleading evidence, is therefore mostly determined by the probability that the hypothesis that is not supported by the current LR is nevertheless true. That probability depends on the prior odds and the current LR.
For example, suppose a large LR for $A$ and $B$ to be siblings is obtained based on their profiles, so that the posterior odds are a million to one. If the profiles of $A$'s parents would now be included for direct comparison with $B$, we would expect the LR to change dramatically towards unrelatedness only when $A$ and $B$ were, in fact, unrelated after all. But that is unlikely, in view of the odds obtained with the evidence that has been duly processed.

%We close with a few remarks. First of all, one should choose the LR as corresponding to the actual knowledge. It is not allowed to compute various LRs and report the most discriminatory one.

%Image that an evaluator uses a black box algorithm that is known to produce accurate LRs in the sense of being calibrated. Now imagine next that the evaluator does not know how the algorithm works, but only that it produces calibrated LRs. Does this matter for whether the evaluator should apply the algorithm? With the arguments in this paper, we see that it does not. When the output of this system is viewed as a score, then the LR of the score is the score itself. This score carries all the evidential information that the evidence processed by the algorithm, whatever it is, carries. 
%However, if may also be that the evaluator learns which data the algorithm uses, and that this is not all data. Should the evaluator now dismiss the algorithm? The output of the algorithm is still a likelihood ratio. Also, if the remaining data would be evaluated, they would be unlikely to yield a large change in LR in the direction of the false hypothesis: a change by a factor $t$ or more has probability at most $1/t.$ The algorithm remains just as applicable, provided we keep ourselves uninformed about the data that it ignores.

So far, we have argued that processing less data than the full data, by using score-based and/or common source LRs, is not in itself a problematic procedure in the sense that it does lead to a decrease in costs of Bayes decisions. That is to say that if an evaluator has only these data, they are better than nothing.
But of course, it is often the case that all data, including the part that was not evaluated, have been observed by the evaluator.  Then additional attention is warranted. Sometimes, also without a quantitative statistical model for the unevaluated data, it may nevertheless be possible to make qualitative statements and to recognize the rare cases where the extra information does have a strong impact.

An example of such a situation has long existed for the interpretation of DNA mixture profiles: such a profile consists of a set of peaks whose location reveals the genomic variants (called the alleles) and whose height is a measure for the abundance of that allele post PCR and hence also pre PCR. Initially, models that gave a LR based on the recorded alleles were available, but these models could not process peak height data. Thus, instead of the full data, only a part (i.e., the observed alleles) was processed in the resulting LR, while the forensic analysts and interpretators had access to the whole profile. They would then, visually and qualitatively and not leading to additional quantitative LR assessment, inspect the correspondence between a person of interest (PoI) and the trace profile. This was done to assess whether they would estimate, according to their qualitative knowledge and expertise, that including the peak heights would provide further  support for contribution of the PoI, if a LR in favour of that hypothesis had been obtained based on the observed alleles without peak heights. Or, conversely, they would first by visual inspection taking alleles and peak heights into account, have to be convinced of the existence of such support, before proceeding to a computation with the statistical model that ignored the peak heights. Of course, if it would be discovered that the peak heights were inconsistent with contribution, while the computed LR supported it, the computed LR would be overruled. This way,  the full data were used, but qualitatively only.

It is of course easy to conceive of peak heights that are inconsistent with the PoI's contribution (and thus would reduce the LR to zero) even when the profile of a PoI and observed alleles of a trace profile lead to a large LR.  However, this is, as in the examples for siblings before, mitigated in practice by the fact that when strong evidence is obtained, this usually means that the posterior probability of contribution is large, and for actual contributors the changes in LRs towards non-contribution are bounded by \eqref{H1bound}. Thus, ignoring the peak heights has primarily the effect of losing evidential strength in the direction of the true hypothesis. This is also why, when a peak-height based model is finally introduced, there is no need to re-visit old cases where strong evidence had been obtained. Only in case where the evidence was weak, this may be worthwhile, time and resources permitting.

Thus, we do not believe that any LR model should be discredited on the basis of the information that it does or does not process. It may in fact be preferable to have a more accurate model for a simplification of the data, than a flawed model for the full data. This can perfectly go together with the evaluator also not being discharged of trying to interpret all relevant data they have access to, even if only qualitatively, and to be convinced of the applicability of the statistical model. This, however, is true for any method, evidence type, and hypotheses.

\section*{Appendix}
\label{sec_proper}

Our analysis has not made use of proper scoring rules, but for completeness we discuss some of the relevant aspects in this appendix.
In \cite{vergeer}, Vergeer ranked various ``LR systems'' for two hypotheses including score based methods, the common source scenario and the specific source scenario. He also argues against the claim of \cite{neumann} that certain methods should not be used, albeit in a rather indirect way using (strictly) proper scoring rules ((S)PSRs) \cite{degroot, dawid2007, brummer2006, raftery, brummer, ferrer}. Using PSRs, it is also possible to obtain the results that we derived directly in section \ref{sec_general}. Although strictly speaking not necessary for this paper, we briefly discuss SPRSs for completeness.

A PSR is a function $C(H_i,q)$ where $q$ is a probability distribution over the $H_i$. The quantity $C(H_i,q)$ is to be interpreted as a cost that applies when we state probability distribution $q$, and then $H_i$ is revealed to be true. The cost $C(H_i,q)$ need not depend on the whole distribution $q$; a well known example is \begin{equation}
\label{logexample}
C(H_i,q)=-\log(q_i),
\end{equation}
which penalizes the occurrence of a hypothesis which was deemed to have small probability.

Given a scoring rule and two probability distributions $q, q'$, we can evaluate the sum (with a slight abuse of notation also denoted by $C$)
\begin{equation}
	\label{fqq}
	C(q' \mid q):=\sum_i q_i C(H_i, q').
\end{equation} 
	
This sum can be interpreted as the expected cost that we incur when we state (or believe) probability distribution $q'$ for the $H_i$, but then the events $H_i$ materialize according to distribution $q$. We call the scoring function (strictly) proper when, for any fixed $q$, $C(q'\mid q)$ is minimal (only) when $q'=q$. In particular, a forecaster tasked with giving a probability distribution for the $H_i$, whilst evaluated by strictly proper scoring rule $C$, will in their own expectation perform best when they offer their own subjective assessment $q$. In general, scoring rules do not have that property: for example, if only the materialization of $H_1$ comes with costs (namely, when the forecaster did not predict it with certainty), then a forecaster wishing to avoid costs will always predict $H_1$ with certainty irrespective of their actual belief. 

SPSRs occur naturally in the context of costs $c_{ij}$, by considering the probability distribution on the $H_i$ that the Bayes decision is based on. 
Indeed, if we assume that we always make a BD, a distribution $q$ induces expected costs
\begin{equation}\label{verwachting}
	\sum_{j=1}^n c_{i_qj}q_j,
\end{equation}
where $i_q$ represents the argmin of \eqref{verwachting}. In this way, taking BDs allows us to define a scoring rule $C_{BD}$ by
\begin{equation}
	\label{spsrbd} C_{BD}(H_j, q)=c_{i_q j}
\end{equation} 
The function $C$ is a proper scoring rule, since
\begin{equation}
	C_{BD}(q \mid q) = \sum_{j=1}^n c_{i_qj}q_j \leq \sum_{j=1}^n c_{i_pj}q_j = C_{BD}(p\mid q).
\end{equation}

In fact, this construction conceptually amounts to considering not only the `hard decisions' $A_i$, but also the `soft decision' preceding it, namely the formulation of the updated probability distribution. Whether we regard the costs of hard decisions (measured by the $c_{ij}$) or the cost of the probability distributions obtained (measured by the associated PSR $C_{BD}$) is immaterial, we measure the same quantity.
Indeed, for any $p$ we have
\[ c(p)=C_{BD}(p\mid p),\]
both representing the expected costs of a Bayes Decision based on $p$, when indeed the events $H_i$ materialize according to $p$.

The inequality in Theorem \ref{bdimp} therefore can be phrased as, for $C=C_{BD}$ derived from costs $c_{ij}$,
\begin{equation}
	\label{mainspsr}
%\BE_e[C((\pi \mid e) \mid (\pi \mid e))] \leq C(\pi).
E_e(C(\pi(e) \mid \pi(e))) \leq C(\pi).
\end{equation}
Therefore, a proof of \eqref{mainspsr} valid for any SPSR $C$ will imply Theorem \ref{bdimp}. This is the approach taken in \cite{vergeer} and \cite{brummer}. In fact, \cite{vergeer} worked in the special case of Section \ref{sec_first} (two hypotheses, two actions) and then treated several special cases that essentially differ only in the choices of $\pi$ and $e$ and then proceed in analogous ways. In \cite{brummer} the general case is treated, but in a rather technical way relying on scoring rules.

There is an interesting relation between proper scoring rules and entropy.
In information theory, the entropy ${\rm Ent}(\pi)$ of a discrete probability distribution $\pi$ is defined as 
$$
{\rm Ent}(\pi) = -\sum_j\pi_jlog(\pi_j).
$$ 
It has multiple interpretations, one of them as a measure of the amount of uncertainty in $\pi$. With this in mind one expects it to decrease (or stay the same) when new data $e$ are conditioned on. Indeed, applying \eqref{mainspsr} to the strictly proper scoring rule in \eqref{logexample}, we see that
$$
\BE\left[ -\sum_{j=1}^n \pi(E)_j \log (\pi(E)_j)\right] \leq -\sum_{j=1}^n \pi_j \log (\pi_j).
$$

For the case $n=2$, if we sample $e$ and obtain likelihood ratio $x=LR(e)$, the posterior distribution will be $\BP(H_1 \mid e)=\frac{x\pi_1}{x\pi_1+\pi_2}$ and $\BP(H_2 \mid e)=\frac{\pi_2}{x\pi_1+\pi_2}$.

The entropy of the posterior distribution is therefore
\[ -\frac{x\pi_1}{x\pi_1+\pi_2}\log(\frac{x\pi_1}{x\pi_1+\pi_2})-\frac{\pi_2}{x\pi_1+\pi_2}\log(\frac{\pi_2}{x\pi_1+\pi_2}).\]
 
The expected entropy of the posterior distribution is then, writing $\ell_i(x)$ for the probability of obtaining $LR=x$ under $H_i$:
\begin{eqnarray*} \BE[{\rm Ent}(\pi(E))] &=& \sum_e \BP(e) {\rm Ent}(\pi(e)) \\ &=& \sum_x \BP(LR=x){\rm Ent}(\pi(x)) 
	\\ &=& \sum_x (\pi_1 \ell_1(x)+\pi_2 \ell_2(x)) {\rm Ent}(\pi(x)) \\ &=& \sum_x (\pi_1x+\pi_2)\ell_2(x)\left(-\frac{x\pi_1}{x\pi_1+\pi_2}\log(\frac{x\pi_1}{x\pi_1+\pi_2})-\frac{\pi_2}{x\pi_1+\pi_2}\log(\frac{\pi_2}{x\pi_1+\pi_2}) \right)\\
	&=& \sum_x \ell_2(x)\left( -x\pi_1 \log(\frac{x\pi_1}{x\pi_1+\pi_2}) - \pi_2 \log(\frac{\pi_2}{x\pi_1+\pi_2}) \right)\\
	&=& -\sum_x \pi_1 \ell_1(x) \log(\frac{x\pi_1}{x\pi_1+\pi_2}) -\sum_x \pi_2 \ell_2(x) \log(\frac{\pi_2}{x\pi_1+\pi_2})\\
	&=& -\pi_1 \BE_{H_1} \log(\BP(H_1\mid E))-\pi_2 \BE_{H_2}\log(\BP(H_2\mid E)).
	\end{eqnarray*}
	
For $\pi_1=\pi_2=1/2$, this quantity is called the CLLR in the literature \cite{brummer2006, vanlierop2024}. We see that this is simply the expected entropy of the posterior distribution starting with uniform priors.

\bibliographystyle{authordate1}

\bibliography{References}

\end{document}